\documentclass[a4paper, BCOR=0.0mm, DIV=calc]{scrartcl}
\usepackage[T1]{fontenc}
\usepackage[utf8]{inputenc}
\usepackage[ngerman]{babel}
\usepackage[osf,sc]{mathpazo}
\usepackage{textcomp}
\usepackage{microtype}
\usepackage{graphicx, color}
\usepackage{amsmath, amssymb, amsfonts}
\usepackage{booktabs}
\usepackage{nicefrac}
\usepackage{wasysym,pifont}
\usepackage{tikz}
\usepackage{graphicx}
\usepackage{float}
\usepackage{wrapfig}
\usetikzlibrary{calc,intersections}
\KOMAoptions{twoside=false, twocolumn=false, headinclude=false, footinclude=false, mpinclude=false, pagesize=auto}
\recalctypearea

\setkomafont{section}{\mdseries\scshape\Large}
\addtokomafont{subsection}{\mdseries\scshape}
\setkomafont{title}{\mdseries\scshape}

\begin{document}
\title{Explicit formulas for orthogonal polynomials derived from their difference equation}
\author{Alexander Aycock}
\date{ }
\maketitle

\section*{Abstract}

We solve the difference equation with linear coefficients by the Momentenansatz to obtain explicit formulas for orthogonal polynomials.

\section*{Keywords}

Momentenansatz; Difference equation; Differential equations of infinite order; orthogonal polynomials; hypergeometric series

\section*{Mathematics subjects classification}

33C05; 33C45; 34A35; 39A06; 

\section*{Introduction}

\paragraph*{§1}

Since the foundational works of \textsc{Euler} \cite{E710} and \textsc{Gauß} \cite{Gauss} the importance of the hypergeometric series

\[
F(a,b,c,z)= \sum_{n=0}^{\infty} \frac{(a)_n(b)_n}{(c)_n}\frac{z^n}{n!}
\]
with

\[
(a)_0=1 \quad \text{and} \quad(a)_n=a(a+1) \cdots (a+n-1)
\]
is well-known, because for example nearly all functions important for physics are special cases of it. 

\paragraph*{§2}

Since \textsc{Euler} it is also known, that the hypergeometric series satisfies the following differential equation

\[
z(1-z) \frac{d^2F}{dz^2}+[c-(a+b+1)z]\frac{dF}{dz}-abF=0
\]
And since \textsc{Gauss'} work on the same it is familiar, that $F$ can be expressed a difference equation of second order with linear coefficients. So in the mentioned paper among many other formulas one finds the following:

\[
0=(c-2a-(b-a)z)F(a,b,c,z)+a(1-z)F(a+1,b,c,z)-(c-a)F(a-1,b,c,z)
\]
which, as \textsc{Gauss} calls it, contiguous relation can be interpreted as difference equation of second order in the variable $a$. If one solves this same equation, one would obtain an explicit representation for $F$.
\paragraph*{§3}

And this paper is dedicated to the resolution of similar difference equations. For this we will use the Momentenansatz, tracing back to \textsc{Euler} (see \cite{E123} or \cite{E594}), which name comes from the similarity to the moment problem, solved by \textsc{Stieltjes}. See \cite{Stieltjes}. This method will allow the resolution of the following difference equation

\[
(\alpha x+a)f(x+2)=(\beta x+b)f(x+1)+(\gamma x +c)f(x)
\]
which contains the example of the hypergepmetric series as a special case, of course. At first we want to explain this method in general and illustrate it with well-known examples.

\section*{First section: Explaination of the Momentenansatz}

\paragraph*{§4} It is propounded to us, to solve the following equation

\[
(\alpha x+a)f(x+2)=(\beta x+b)f(x+1)+(\gamma x +c)f(x)
\]
with arbitrary numbers $\alpha, a, \beta, b, \gamma, c$. Or completely general

\[
(\alpha_n x+a_n)f(x+n)=\sum_{k=0}^{n-1}(\alpha_k x+a_k)f(x+k)
\]
where we want to focus mainly on the case $n=2$. The generalisation to the case $n>2$ is neither difficult nor neccessary for our purpose. So all the following will hold only for the case $n=2$.

\subsection*{Remark}

\paragraph*{§5} Depending on the nature of the variable $x$ the number of the constants and therefore the neccessary initial conditions for an unique solution varies. If $x \in \mathbb{Z}$, so we need two initial conditions.\\
If $x \in \mathbb{C} \setminus \mathbb{Z}$, so one needs infintely many conditions. Hence we have for holomorphic $f$

\[
f(x+a)= \sum_{n=0}^{\infty} \frac{f^{(n)}(x)}{n!}a^n
\]
so especially

\begin{alignat*}{8}
&f(x+1)&&=\sum_{n=0}^{\infty} \frac{f^{(n)}(x)}{n!}\\
&f(x+2)&&=\sum_{n=0}^{\infty} \frac{f^{(n)}(x)}{n!}2^n
\end{alignat*}
Therefore we can see

\[
(\alpha x +a)f(x+2)=(\beta x +b)f(x+1)+(\gamma x +c)f(x)
\]
as differential equation of infinite order.\\
From this we see, that number of conditions is indeed infinite, but countable. Furthermore we recongnize, that the equation for $x \in \mathbb{Z}$ is interpolated this way. This will be important later. But for now we start with the definition of the Momentenansatz.

\subsection*{Definition (Momentenansatz)}

\paragraph*{§6} Let a equation of the form from above be propounded, then we call the ansatz

\[
f(x)= \int\limits_{a}^{b} R(t)^{x-1}h(t)dt
\]
a Momentenansatz. Here $R(t)$ and $h(t)$ are two arbitrary functions, which will later be determined from the concrete problem.

\subsection*{Remark}

\paragraph*{§7} From the Momentenansatz it follows, that we have to determine the following quantities

\begin{alignat*}{9}
&\text{I.} \quad  && \text{The function} \quad && R(t)\\
&\text{II.} \quad  && \text{The function} \quad && h(t)\\
&\text{III.} \quad  && \text{The boundaries} \quad &&a ~ \text{and} ~ b\\
\end{alignat*}
For this we introduce the following auxiliary equation:

\subsection*{Definition (Auxiliary Equation)}

\paragraph*{§8} If the equation

\[
(\alpha x +a)f(x+2)=(\beta x +b)f(x+1)+(\gamma x +c)f(x)
\]
was propunded, we call the equation

\[
(\alpha x+a) \int\limits_{}^{t} R^{x+1}(t)h(t)dt=(\beta x+b) \int\limits_{}^{t} R^{x}(t)h(t)dt+(\gamma x+c) \int\limits_{}^{t} R^{x-1}(t)h(t)dt+R(t)^xQ(t)
\]
the auxiliary equation. 

\[
F(t)= \int\limits_{}^{t} f(t)dt
\]
is the integral function of $f$. $Q(t)$ is another function to be determined.

\subsection*{Remark}

\paragraph*{§9} Although it might seem absurd, to introduce another function $Q(t)$, we will see in one moment, that exactly by this we will be able to solve the equation.\\
Hence it is possible, to derive a system of differential equations from the auxiliary equation, that relates the functions $h(t)$, $Q(t)$ and $R(t)$. Thereafter it will also be possible, to determine the boundaries.

\subsection*{Definition (Resolvent)}

\paragraph*{§10} We call conditions following from the auxiliary equation the Resolvent of the Momentenansatz. In the next paragraphs we want to derive these conditions. 

\subsection*{Derivation of the conditions for $Q(t)$, $R(t)$ and $h(t)$}

\paragraph*{§11}

Let us start from our auxiliary equation

\[
(\alpha x+a) \int\limits_{}^{t} R^{x+1}(t)h(t)dt=(\beta x+b) \int\limits_{}^{t} R^{x}(t)h(t)dt+(\gamma x+c) \int\limits_{}^{t} R^{x-1}(t)h(t)dt+R(t)^xQ(t)
\]
If we differentiate this equation with respect to $t$ and use the Fundamental theorem of Calculus, we obtain

\[
(\alpha x+a) R^{x+1}(t)h(t)= (\beta x +b)R(t)^{x}h(t)+(\gamma x+c)  R^{x-1}(t)h(t)+xR^{x-1}(t)Q(t)R^{\prime}(t)+R^x(t) Q^{\prime}(t)
\]
And after the division by $R(t)^{x-1}$

\[
(\alpha x+a) R^2(t)h(t)= (\beta x +b)R(t)h(t)+(\gamma x+c)  h(t)+xQ(t)R^{\prime}(t)+R(t) Q^{\prime}(t)
\]
If we compare  the coefficients of $x$ on both sides of the equation, we find the following system of equations:

\begin{alignat*}{9}
& \text{I.} \quad && \alpha R^2(t)h(t)&&=\beta R(t) h(t) &&+\gamma h(t) &&+R^{\prime}(t)Q(t) \\
& \text{II.} \quad && a R^2(t)h(t)&&=b R(t) h(t) &&+c h(t) &&+Q^{\prime}(t)R(t) \\
\end{alignat*}
which is a system of differential equations for the functions $h(t)$, $Q(t)$ and $R(t)$.

\subsection*{Corollary 1}

\paragraph*{§12} The system of equations is under-determined, because one has only two equations for three functions. Hence one take one of the theree functions completely ad libitum. One has the three following possibilities:

\begin{alignat*}{9}
&\text{I.} \quad && h(t) \quad && \text{is arbitrary, hence}
\quad R(t) &&\quad \text{and} \quad Q(t) && \quad \text{are functions of} && \quad h(t)  \\
&\text{II.} \quad && Q(t) \quad && \text{is arbitrary, hence}
\quad R(t) &&\quad \text{and} \quad h(t) && \quad \text{are functions of} && \quad Q(t)  \\
&\text{III.} \quad && R(t) \quad && \text{is arbitrary, hence}
\quad h(t) &&\quad \text{and} \quad Q(t) && \quad \text{are functions of} && \quad R(t)  \\
\end{alignat*}
Which function is put ad libitum, follows from the particular problem in the most cases.

\subsection*{Corollary 2}

\paragraph*{§13} If one chooses $h(t)$ ad libitum, so one has two arbitrary constants, because $R(t)$ and $Q(t)$ appear in their first derivative.\\
If one on the other hand chooses one of the other functions ad libitum, only one constant appears.

\subsection*{Remark}

\paragraph*{§14}

From the last corollary we see, that for the case $x \in \mathbb{Z}$, in which two initial conditions have to be given for the uniqueness of the problem, we can proceed in two ways. We can either take $h(t)$ ad libitum and solve the problem directly. Or we choose $R(t)$ or $Q(t)$ ad libitum, with that restriction, that the second condition is met automatically.\\
The case $x \in \mathbb{C} \setminus \mathbb{Z}$ is a lot harder and one only obtains a special solution, if it is not possible, to satisfy all initial conditions at once. \\
Finally we want to explain how we can determine the integral boundaries.

\subsection*{Definition (Determinant)}

\paragraph*{§15}

We call the equation

\[
R^x(t)Q(t)=0
\]
the Determinant for the boundries of the integral

\[
\int\limits_{}^{t} R^{x-1}(t)h(t)dt
\]

\subsection*{Corollary}

\paragraph*{§16}

So if the Determinant has less than $2$ solutions, so obtains no solutions with the Momentenansatz.\\
If it has exactly $2$ solutions, so one finds exactly one solution.\\
But if it are $n>2$ solutions, so one has $\binom{n}{2}$ solutions of the problem.

\subsection*{Remark 1}

\paragraph*{§17}

The definition of the Determinant immediately follows from the auxiliary equation

\[
(\alpha x+a) \int\limits_{}^{t} R^{x+1}(t)h(t)dt=(\beta x+b) \int\limits_{}^{t} R^{x}(t)h(t)dt+(\gamma x+c) \int\limits_{}^{t} R^{x-1}(t)h(t)dt+R(t)^xQ(t)
\]
For, if one finds two values $a$ and $b$ of $t$, so that the Determinant

\[
R^x(t)Q(t)=0
\]
is satisfied, so one finally finds

\[
(\alpha x+a) \int\limits_{a}^{b} R^{x+1}(t)h(t)dt=(\beta x+b) \int\limits_{a}^{b} R^{x}(t)h(t)dt+(\gamma x+c) \int\limits_{a}^{b} R^{x-1}(t)h(t)dt+0
\]
and therefore

\[
f(x)=\int\limits_{a}^{b} R^{x-1}(t)h(t)dt
\]
as solution of the equation

\[
(\alpha x +a)f(x+2)=(\beta x+b)f(x+1)+(\gamma x+c)f(x)
\]
hence of the propounded problem itself.

\subsection*{Remark 2}

\paragraph*{§18}

In certain problems it can be convenient, to put instead of

\[
f(x)= \int\limits_{}^{t} R^{x-1}(t)h(t)dt
\]
more general

\[
f(x)= \int\limits_{}^{t} R^{x+k}(t)h(t)dt,
\]
what becomes clear in the peculiar cases and only shortens the calculations, but does not alter the procedure in the solution of the problem. The steps are always the following:\\
1. Make the Momentenansatz\\
2. Formulate the auxiliary equation\\
3. Derive the Determinant\\
4. Choose one of the 3 functions ad libitum an solve the system of equations\\
5. Solve the Determinant\\
6. Determine the constants from the initial condition\\
To 4. we remark, that the function should be taken in such a way, that can be done, so that the Determinant has a solution.\\
Now having explained everything, we want to illustrate it with examples.

\section*{Second Section: Examples}

\subsection*{Example 1: The Gamma Function}

\paragraph*{§19} Let the following equation be propounded

\[
\Gamma(x+1)=x\Gamma(x) \quad \text{with} \quad \Gamma(1)=1
\]
We illustrate the procedure with all steps:\\
$\mathbf{Step ~ 1}$: The Momentenansatz\\
We put, that we have

\[
\Gamma(x)= \int\limits_{a}^{b}R^{x-1}(t)h(t)dt
\]
$\mathbf{Step ~ 2}$: Auxiliary Equation\\
We formulate the following equation

\[
\int\limits_{}^{t}R^{x}(t)h(t)dt=x\int\limits_{}^{t}R^{x-1}(t)h(t)dt+R^x(t)Q(t)
\]
$\mathbf{Step~ 3}$: Derivation of the Resolvent\\
We differentiate the auxiliary equation, divide by $R^{x-1}(t)$ and find:

\[
R(t)h(t)=xh(t)+xR^{\prime}(t)Q(t)+R(t)Q^{\prime}(t)
\]
A comparison of the coefficients of the powers of $x$ gives the following system:

\begin{alignat*}{9}
&\text{I.} \quad &&R(t)h(t) &&=R(t)Q^{\prime}(t)\\
& \text{II.} \quad && \quad ~~0 &&=+h(t)+R^{\prime}(t)Q(t)
\end{alignat*}
$\mathbf{Step ~ 4}$: Choosing a function ad libitum and solving the system\\
Hence the system is easily solved, we choose

\[
R(t)=t
\]
The system then becomes

\begin{alignat*}{9}
&\text{I.} \quad &&h(t) &&=Q^{\prime}(t)\\
& \text{II.} \quad &&  h(t) &&=-Q(t)
\end{alignat*}
and is solved by

\[
Q(t)=Ce^{-t} \quad \text{und} \quad h(t)-Ce^{-t}
\]
$\mathbf{Step ~ 5}$: Solution of the Determinant\\
The Determinant is

\[
0=t^xQ(t)=Ct^xe^{-t}
\]
For $x>0$ one finds the two solutions

\[
t=0 \quad \text{and} \quad t= \infty
\]
This leads to the intermediate result

\[
\Gamma(x)= -C\int \limits_{0}^{\infty}e^{-t}t^{x-1}dt
\]
$\mathbf{Step ~ 6}$: Determination of the constant of integration\\
According to the initial condition it has to be

\[
1=\Gamma(1)=-C \int \limits_{0}^{\infty}e^{-t}dt=-C \cdot 1
\]
Hence

\[
C=-1
\]
And our solution then is

\[
\Gamma(x)= \int \limits_{0}^{\infty} e^{-t}t^{x-1}dt
\]
\subsection*{Remark}

\paragraph*{§20}

Hence we just derive the integral representation of $\Gamma(x)$, even the condition for the convergence of the integrals appeared by itself.\\
Further we note, that the choice $R(t)=t$ at the same time gave the most simple soultion of the difference equation

\[
\Gamma(x+1)=x\Gamma(x).
\]
If one considers  $R(t)$ in general, one would find

\[
\Gamma(x)= \int \limits_{0}^{\infty} e^{-\int\limits_{}^{t}R(t)dt}R^{x-1}(t)R^{\prime}(t)dt
\]
where $R(t)$ would just have to satisfy the following conditions

\[
R(0)=0 \quad \text{and} \quad R(\infty)=\infty
\]
That in some way $R(t)=t$ is a convenient choice, will be seen in the following examples.\\
It illustrates the always reappearing phenomenon in Physics, that the most simple solution for a problem with several solutions is always the physically relevant one. Furthermore the Momentenansatz 

\[
\Gamma(x)= \int \limits_{a}^{b} h(t)t^{x-1}dt
\]
seems to be a promising ansatz also for other physical problems.

\subsection*{Example 2: Die Legendre-Polynomials}

\paragraph*{§21} The Legendre-Polynomials satisfy the following difference equation

\[
P_{n+1}(x)=(2n+1)xP_n(x)-nP_{n-1}(x)
\]
Note, that the difference equation is one for $n$, not $x$. The first two Legendre-Polynomials are $P_0=1$ and $P_1=x$.\\
Encouraged by the physical principle, we put:

\[
P_n(x)=\int\limits_{a}^{b} t^n h(t)dt
\]
so put $R(t)=t$ directly. The auxiliary is 

\[
\int\limits_{}^{t} t^{n+1} h(t)dt=(2n+1)x\int\limits_{}^{t} t^n h(t)dt-n\int\limits_{}^{t} t^{n-1} h(t)dt+t^nQ(t)
\]
and leads to the Resolvent

\begin{alignat*}{9}
&\text{I.} \quad &&h(t)(t^2-2xt+1)&&=Q(t)\\
&\text{II:} \quad &&h(t)(t^2-xt)&&=tQ^{\prime}(t)
\end{alignat*}
with the solutions

\[
Q(t)=C\sqrt{t^2-2xt+1} \quad \text{und} \quad h(t)= \frac{C}{\sqrt{t^2-2xt+1}}
\]
So the Determinant is

\[
0=Q(t)t^n=Ct^n\sqrt{t^2-2xt+1}
\]
with the solutions

\[
t=x+\sqrt{x^2-1} \quad \text{und} \quad t=x-\sqrt{x^2-1}
\]
This yields the intermediate solution

\[
P_n(x)=C\int\limits_{x-\sqrt{x^2-1}}^{x+\sqrt{x^2-1}} \frac{t^n}{\sqrt{1-2xt+t^2}}dt
\]
and because of $P_0(x)=1$ we find:

\[
1=C\int\limits_{x-\sqrt{x^2-1}}^{x+\sqrt{x^2-1}} \frac{dt}{\sqrt{1-2xt+t^2}}=C\bigg[\log(\sqrt{t^2-2xt+1}+t-x)\bigg]_{x-\sqrt{x^2-1}}^{x+\sqrt{x^2-1}}=C\log(-1)=Ci \pi
\]
Hence

\[
C= \frac{1}{i \pi}
\]
and our general formula.

\[
P_n(x)=\frac{1}{i \pi}\int\limits_{x-\sqrt{x^2-1}}^{x+\sqrt{x^2-1}} \frac{t^n}{\sqrt{1-2xt+t^2}}dt.
\]

\subsection*{Remark 1}

\paragraph*{§22}

It was not surprising, that the variable $x$ appeared in the boundaries, because the difference equations is for $n$, so $x$ is just a constant. It could indeed be doubted, if $P_1(x)=x$, because the recurrence equations for $P_n(x)$ depends on two initial conditions and we only used $P_0(x)=0$ explicitly. But one easily calculates, that our formula yields the correct value. By our ansatz we satisfied one condition automatically. Just note, that we could have taken evrey Legendre-Polynomial as the other initial condition and would have gotetn the same formual then, which is indeed quite remarkable and was not to be exspected.

\subsection*{Remark 2}

\paragraph*{§23}

The formula

\[
P_n(x)=\frac{1}{i \pi}\int\limits_{x-\sqrt{x^2-1}}^{x+\sqrt{x^2-1}} \frac{t^n}{\sqrt{1-2xt+t^2}}dt
\]
also directly allows a definition of $P_n(x)$ for $n \notin \mathbb{N}$; and it would be worth the effort to check, if this then agrees with the usual definitions in the literature. Note, that all values $\frac{n}{2}$ with $n \in \mathbb{Z}$ leads to elliptic integrals. But we do not want to consider this in more detail here.

\subsection*{Example 3: The Hermite-Polynomials}

\paragraph*{§24}

The Hermite-Polynomials satisfy the following equation

\[
H_{n+1}(x)=2xH_n(x)-2n H_{n-1}(x) \quad \text{with} \quad H_=(x)=1, \quad H_1(x)=2x
\]
Ansatz and auxiliary equation are determined in the same way as above and lead to the following Resolvent

\begin{alignat*}{9}
&\text{I.} \quad &&=\frac{1}{2}Q(t)\\
& \text{II.} \quad &&=\frac{tQ^{\prime}(t)}{t^2-2xt}
\end{alignat*}
with solution

\[
Q(t)=Ce^{\frac{1}{4}(t^2-4xt)} \quad \text{und} \quad h(t)=\frac{1}{2}Ce^{\frac{1}{4}(t^2-4xt)}
\]
Hence the Determinant is

\[
0=Ct^ne^{\frac{1}{4}(t^2-4xt)}
\]
which is solved by

\[
t=-i \infty \quad \text{und} \quad t=i \infty 
\]
So we have the following intermediate solution

\begin{alignat*}{9}
&H_n(x)&&=C\int\limits_{-i \infty}^{i  \infty} t^ne^{\frac{t^2}{4}-xt}dt \quad \text{and for} \quad t=iy\\
& &&=\frac{C}{i^{n-1}}\int\limits_{- \infty}^{  \infty} t^ne^{\frac{-y^2}{4}-ixy}dy 
\end{alignat*}
And because of initial condition it has to hold

\[
1=H_0(x)=i C\int\limits_{-\infty}^{  \infty}e^{\frac{-y^2}{4}-ixy}dy=2i C \sqrt{\pi}e^{-x^2}
\]
and therefore

\[
C= \frac{e^{x^2}}{2 i \sqrt{\pi}}
\]
So our formula is

\begin{alignat*}{9}
&H_n(x)&&=\frac{e^{x^2}}{2\sqrt{\pi}i^n} \int\limits_{- \infty}^{  \infty} t^ne^{\frac{-y^2}{4}}e^{-ixy}dy\\
& &&=\frac{e^{x^2}}{\sqrt{2}i^n}F(y^ne^{-\frac{y^2}{4}})
\end{alignat*}
where $F$ means the Fourier Transform here. \\
In exactly the same way one finds for the Laguerre-Polynomials with the difference equation

\[
(n+1)L_{n+1}(x)=(2n+1-x)L_n(x)-nL_{n+1}(x) \quad \text{mit} \quad L_0(x)=1, \quad L_1(x)=1-x
\]
the formula

\[
L_n(x)= \frac{1}{Ei(x)}\int_0^1 \frac{t^ne^{\frac{-x}{t-1}}}{t-1}dt
\]
with the exponential integral $Ei(x)$. Here again the phenomenon occurs, that we only needed one initial condition explicitly, to derive the formula. In addition this formula does not seem to appear in the literature. It is of great interest, how this can be explained mathematically and if this can be generalised to other classes of functions. 

\subsection*{Remark 1}

\paragraph*{§25} Now we were able to show the method with some important examples and see quickly, that the difference equation for the hypergeometric series could be treated in the same way. Without great effort one reaches the well-known integral representation

\[
F(a,b,c,z)= \frac{\Gamma(c)}{\Gamma(b)\Gamma(c-b)}\int\limits_{0}^{1}t^{b-1}(1-t)^{c-b-1}(1-tz)^{-a}dt
\]
which was proved by \textsc{Euler}. \cite{E366}.\\
And because the hypergeometric function because of the Sturm-Liouville is of greatest importance, these explainations are maybe not without use.

\paragraph*{§26}

In addition is was remarkable, that the ansatz $R(t)=t$, so the most simple possible, always led to the correct result and contained all physically relevant properties. What is the interpretation of the general ansatz, needs a further investigation. But the follwing principle is almost more important for us.

\subsection*{Priciple of the simplest solution}

If there are more solutions to a mathematical problem, so the most simple non-trivial one is always the one important for physics.

\end{document}